\documentclass[openacc]{rstransa}

\newtheorem{theorem}{\bf Theorem}[section]

\DeclareMathOperator*{\loc}{loc}
\DeclareMathOperator*{\ess}{ess}

\DeclareMathOperator*{\Ker}{Ker}

\DeclareMathOperator*{\R}{Re}
\DeclareMathOperator*{\I}{Im}

\newcommand{\inv}{^{-1}}

\newcommand{\RR}{\mathbb{R}}
\newcommand{\CC}{\mathbb{C}}

\newcommand{\T}{(T(t))_{t\ge 0}}

\newcommand{\Mlog}{M_{\mathrm{log}}}

\theoremstyle{definition}

\begin{document}

\title{Semi-uniform stability of operator semigroups and energy decay of damped waves}

\author{
R. Chill$^{1}$, D. Seifert$^{2}$ and Y. Tomilov$^{3}$}

\address{$^{1}$Institut f\"{u}r Analysis, Fakult\"{a}t f\"{u}r Mathematik, TU Dresden, 01062 Dresden, Germany
\email{ralph.chill@tu-dresden.de}\\
$^{2}$School of Mathematics, Statistics and Physics, Newcastle University, Newcastle upon Tyne, NE1 7RU United Kingdom
\email{david.seifert@ncl.ac.uk}\\
$^{3}$Institute of Mathematics, Polish Academy of Sciences, \'{S}niadeckich 8, 00956 Warsaw, Poland
\email{ytomilov@impan.pl}}

\subject{xxxxx, xxxxx, xxxx}


\corres{R. Chill\\
\email{ralph.chill@tu-dresden.de}}

\begin{abstract}
Only in the last fifteen years or so has the notion of semi-uniform stability, which lies between exponential stability and strong stability, become part of the asymptotic theory of $C_0$-semigroups. It now lies at the very heart of modern semigroup theory. After briefly reviewing the notions of exponential and strong stability, we present an overview of some of the best known (and often optimal) abstract results on semi-uniform stability. We go on to indicate briefly how these results can be applied to obtain (sometimes optimal) rates of energy decay for certain damped second-order Cauchy problems.
\end{abstract}


\begin{fmtext}
\section{Introduction}\label{sec:intro}

Exponential and strong stability of $C_0$-semigroups are two classical topics in semigroup theory, and the literature on these topics,
through various deep results over the past fifty years, has now reached a reasonably complete state; we refer to \cite{Na86,Ne96,ABHN01} for extensive accounts. Exponential stability is a  strong property, and
it has a number of natural applications arising from  its specific quantitative character and its robustness under perturbations. Meanwhile strong stability, that is to say mere convergence  to zero of semigroup orbits, is a rather delicate property which (in the absence of exponential stability) is highly sensitive to perturbations and depends on fine properties of the semigroup generator.

\end{fmtext}


\maketitle

The present survey aims to address recent progress on how the gap between these two extreme kinds of semigroup stability can be bridged.
After reviewing some basic aspects of exponential and strong stability and, where appropriate, offering some of our own commentary,
we describe an abstract approach to the study of quantified stability for operator semigroups based on Tauberian theorems for Laplace transforms with remainder terms.
This allows us to view strong stability as an end-point case of quantified stability theory,
with exponential stability at the opposite end.
Even though the pioneering work in this area goes back at least fifteen years, certain aspects of the theory may not be well known even among experts.  Our main message is simple: many types of quantified asymptotic behaviour of operator semigroups can be characterised, or at least described very precisely, 
in  terms of simple resolvent bounds for semigroup generators and of various related resolvent conditions. 
Strikingly, the passage from resolvent bounds to decay rates can often be achieved without essential loss, so that optimal decay rates are possible as long as one has sharp resolvent estimates. In fact,   sharp resolvent estimates are (almost) equivalent to sharp decay rates. While resolvents have long been used in the study of decay rates,
a unified approach leading to optimal rates of decay has emerged only relatively recently.  The techniques outlined in the present survey are intended to serve as a partial remedy,
although they too have their natural limitations and will not be the optimal tool in all situations.

We begin, in Section~\ref{sec:exp_str}, by revisiting the classical subjects of exponential and strong stability, paying particular attention to their connections with other topics in the asymptotic theory of operator semigroups such as spectral mapping theorems. Then, in Section~\ref{sec:semi}, we turn to the modern theory of semi-uniform stability of operator semigroups, approached through quantified Tauberian theory before finally, in Section~\ref{damped}, outlining how the abstract theory can be applied to the study of abstract second order problems and thus to energy decay of damped waves. There are many other interesting and important applications of the theoretical results, for instance to decay of correlations of chaotic flows on manifolds, which due to space limitations we are unable to address in our survey. Moreover, even the selection of theoretical topics covered here is necessarily incomplete and rather skewed towards optimal results on Hilbert spaces, which are central to the study of damped waves. Nevertheless, we hope that the survey will shed new light on both classical and modern aspects of the quantified asymptotic theory of $C_0$-semigroups.

\section{Exponential and strong stability revisited}\label{sec:exp_str}
\subsection{Exponential stability and spectral mapping theorems}\label{expon_stab}

A $C_0$-semigroup $(T(t))_{t \ge 0}$ on a Banach space $X$ is \emph{exponentially stable} if for every $x \in X$ there exist $M(x)>0$ and $\omega(x)>0$
such that $\|T(t)x\| \le M(x)e^{-\omega(x) t}$ for all $t \ge 0$. By a simple argument based on the uniform boundedness principle, this is equivalent to the existence of positive constants $M$ and $\omega$ not depending on $x$ such that $\|T(t)\|\le Me^{-\omega t}$, $t \ge 0$. Thus if one defines the \emph{exponential growth bound} of $\T$ as 
\begin{equation}\label{omega}
 \omega (T) := \inf \big\{ \omega\in\RR : \exists M\geq 0\mbox{ s.t.\ } \| T(t)\|\leq M e^{\omega t}\mbox{ for all }t\ge0\big\}, 
\end{equation}
then exponential stability is equivalent to the property $\omega(T) <0$.
Of course, if there exists $t_0>0$ such that $\|T(t_0)\|< 1$, then $\T$ is exponentially stable (and vice versa).
Thus, any uniform decay rate of the form $\|T(t)x\| = O(r(t))$ for some $r\in C_0 (\RR_+)$ and every $x \in X$  implies exponential stability of $\T$
in view of the uniform boundedness principle. Once $\T$ is known to be exponentially stable one is often interested in the optimal choices of $M$ and $\omega$ above, and we will address this matter in some detail below.

From the point of view of many applications exponential stability is the ``best'' type of stability one may hope for. It is robust under small bounded perturbations, indeed even under certain non-linear perturbations, and as a result it forms the basis for the theory of linearised stability. Equilibrium points with exponentially stable linearisations  are easily detected, whereas equilibrium points whose linearisations are stable merely in a weaker sense are invisible in practice; see however \cite{JoLa13,JoLa20}, which deal with an intermediate situation in the spirit of our survey. In many situations the primary task is to characterise exponentially stable or strongly stable semigroups in terms of appropriate a priori information, since the semigroup is rarely known explicitly, and sometimes even when it is the computations quickly become cumbersome. Natural a priori objects, both from the point view of abstract operator theory and also for applications to PDEs, are the spectrum  $\sigma(A)$ of the generator $A$ of a $C_0$-semigroup $(T(t))_{t \ge 0}$ and the resolvent $R(z, A):=(z-A)^{-1}$ for $z$ in the resolvent set $\rho(A):=\CC\setminus\sigma(A)$ of $A$. If $A$ is bounded, then by the elementary spectral mapping theorem for the Riesz-Dunford functional calculus one has
$\sigma(T(t))=e^{t\sigma(A)}$, $t \ge 0$, and the spectral radius formula implies that 
\begin{equation}\label{growth_det}
\omega(T)=s(A),
\end{equation} 
where $s(A) := \sup \{ \R z : z\in\sigma (A) \}$
 is the \emph{spectral bound} of $A$.
In particular, $s(A)<0$ is equivalent to $\omega(T)<0$, and the spectral bound of $A$ alone determines whether or not the semigroup $\T$ is exponentially stable.
Unfortunately, \eqref{growth_det} may fail dramatically if $A$ is an unbounded operator. Recall the fundamental relation from semigroup theory,
\begin{equation}\label{fund}
R(z,A)x=\int_0^\infty e^{-z t}T(t)x\, dt,\quad x\in X,\ \R z  >\omega(T),
\end{equation}
which says that on suitable half-planes the resolvent of the generator is the Laplace transform of the semigroup (in the strong sense).
In particular, the right half-plane $\{z\in\CC:\R z >\omega (T)\}$ is contained in the resolvent set $\rho(A)$, and moreover $s(A) \le \omega(T)$. 
However, one of the main differences between the theory of Laplace transforms and the theory of power series is
that Laplace transforms do not allow for an analogue of the Cauchy-Hadamard formula for power series. Thus on any right half-plane $\mathbb C_\alpha:=\{z\in\CC: \R z \ge \alpha\}$ the integral on the right-hand side of \eqref{fund} may be convergent without being absolutely convergent, and its convergence on $\mathbb C_\alpha$ does not guarantee boundedness of $R(z, A)$ on any smaller half-plane; for relevant examples see e.g.\ \cite{Ba03} and the references therein.

It is natural therefore either to try to describe those situations in which \eqref{growth_det} \emph{does} hold or to describe $\omega(T)$ in non-spectral terms, and there has been extensive activity in both directions.
An obvious place to look for semigroups satisfying \eqref{growth_det} is  among semigroups $\T$ such that  
\begin{equation*}\label{smt}
\text{either} \quad \sigma (T(t))\setminus \{0\}=e^{t\sigma (A)},  t \ge 0, \quad  \text{or at least} \quad \sigma(T(t))=\overline {e^{t\sigma (A)}},  t \ge 0.
\end{equation*}
In the first case, one says that the \emph{spectral mapping theorem} (SMT) holds for $\T$, while in the second case, $\T$ is said to satisfy the \emph{weak spectral mapping theorem} (WSMT).
It has long been known that the SMT holds for {eventually norm-continuous} semigroups (that is to say, $C_0$-semigroups which are continuous in the uniform operator topology on $(t_0,\infty)$ for some $t_0\ge0$); see e.g.\ \cite[Corollary IV.3.12]{EnNa99}, \cite[Theorem~16.4.1]{HiPh57}, \cite[Theorem~2.3.2]{Ne96}. This class is rather large and, in particular, it includes all eventually differentiable semigroups (thus also all analytic semigroups) and all eventually compact semigroups.
There are partial generalisations of the SMT based on replacing eventual norm-continuity by \emph{norm continuity at infinity} and related notions; see e.g.\ \cite{Bl01}, \cite[Theorem~2.3.3]{Ne96}. (Note that on Hilbert spaces norm continuity of a  $C_0$-semigroup on $(0,\infty)$  can be characterised by the simple resolvent condition $\|R(z, A)\|\to 0$ along an appropriate vertical line.)
Unfortunately, these topics are beyond the scope of this survey.
For interesting applications of the spectral mapping theorems to the study of concrete linear and non-linear PDEs
see e.g.\ the recent works \cite{DoSc14,DoSc15,GuMiMo17} and also the slightly less recent papers \cite{CrLa03,GJLS00}.

The WSMT holds for $C_0$-groups $(T(t))_{t\in\RR}$ with sufficiently moderate growth as $|t|\to\infty$  \cite[Theorem~2.4.4]{Ne96}. However, such (semi)groups cannot be exponentially stable. Note also that by a theorem due to Weis \cite{We95} the equality \eqref{growth_det} holds for positive semigroups on $L^p$-spaces and $C(K)$-spaces. 
However, as noted in  \cite{Ta86}, there exists a positive group on the Banach lattice $L^1$ failing to satisfy even the WSMT.
Thus it is possible for the SMT to fail even though the crucial property \eqref{growth_det} still holds. This situation seems to be completely unexplored in the literature. The multiplication semigroup on $\ell^2(\mathbb N)$ given by $T(t)x=(e^{int}x_n)_{n \ge 1}$ for $x=(x_n)_{n \ge 1}$ provides a simple  example of a $C_0$-semigroup satisfying the WSMT but not the SMT. Apart from the one given in \cite{Ta86}, examples of $C_0$-semigroups failing to satisfy the WSMT (with a detailed analysis of the spectrum and resolvent behaviour) may be found in \cite{HeKa86,Wr89}. All of these examples are based on a widely cited example of Zabczyk \cite{Za75} (preceded by \cite{Fo73}) based on a direct sum construction. At the same time, the somewhat neglected example of failure of the WSMT for the Riemann-Liouville $C_0$-group $T((t))_{t \in \mathbb R}$ on $L^2(0,1)$, defined by $(T(t)f)(s)=\frac{1}{\Gamma (t)}\int_{0}^s(s-r)^{it-1} f(r)\, dr$,  can already be found in \cite[Section 23.6]{HiPh57}. We refer to \cite{Hen87} for an illuminating study of this group and of the resolvent of its generator.  Note that in this case one has  $\sigma (T(t))=\{z\in\CC: e^{-|t|\pi/2}\le |z|\le e^{|t|\pi/2}\}$ for $t \in \mathbb R$ while $\sigma(A)=\emptyset$. Moreover, the resolvent of $A$ is compact in this example.

Another line of counterexamples to the SMT stemming from examples due to Wolff and Greiner, Voigt and Wolff is based on the study of translation semigroups like $(T(t)f)(s) = f(e^t s)$ on intersections of weighted $L^p$-spaces or on the Sobolev space $H^1$ over $(1,\infty )$ \cite[Example 5.3.2]{ABHN01}; see  \cite{ABHN01} for a thorough discussion and further references. 
Renardy's paper \cite{Re93} is a standard reference for failure of the WSMT in the case of a semigroup originating from the simple  hyperbolic equation $u_{tt} - u_{xx} - u_{yy} = e^{2\pi xi} u_y$ on $(0,1) \times (0,1) \times \mathbb R$ with periodic boundary conditions; see also \cite{BaLiXi05,RoVe18a}. Other important counterexamples from the point of view of applications are given in Lebeau \cite{Le96}, Schenck \cite{Sc11} and Jin \cite{Jin20} in the context of the damped wave equation $u_{tt} + b u_t - u_{xx} - u_{yy} = 0$.

Recall that the essential spectrum $\sigma_{\ess}(T)$ of a bounded linear operator $T$ can be defined as $\sigma_{{\ess}}(T):=\{z\in\CC: z - T$ is not Fredholm$\}$. It is instructive to observe that the set $\sigma(T)\setminus \sigma_{{\ess}}(T)$ is at most countable and consists of isolated eigenvalues of $T$ of finite multiplicity, accumulating only at points in $\sigma_{{\ess}} (T)$.
Thus $\sigma_{\ess}(T)$ in a natural sense contains ``most'' of the spectrum of $T$. The essential spectrum of Hilbert space $C_0$-semigroups $(T(t))_{t \ge 0}$ originating from a large class of hyperbolic equations was described in non-resolvent, geometric  terms in the deep paper \cite{KoTa95} by Koch and Tataru. More specifically, it was shown in \cite[Theorem~3]{KoTa95} that if $(T(t))_{t \ge 0}$ governs a hyperbolic initial boundary value problem on a smooth, compact and connected Riemannian manifold with smooth boundary, then, under some mild technical assumptions, $\sigma_{\ess}(T(t))$ has rotational symmetry, being either a disk, an annulus or the union of an annulus and the origin for all times $t$ outside a countable exceptional set. In particular, if $(T(t))_{t \ge 0}$ is the semigroup arising in Renardy's example, then $\sigma_{{\ess}}(T(t)) = \{z\in\CC : e^{-t/2} \le |z| \le e^{t/2} \}$, while the spectrum of the generator lies on the imaginary axis; see \cite[Example 2]{KoTa95}.  It would be instructive to develop a more complete theory involving the properties of the resolvent of the semigroup generator. We point out that the spectral description due to Koch and Tartaru is related to the famous ``circles conjecture'' from semigroup theory; see e.g. \cite[Theorem~4.1]{Wr89} and the comments following it. In \cite{He82spec}, Herbst proved that if the SMT fails in the Hilbert space setting then it fails dramatically. Specifically, Herbst showed that if $(T(t))_{t \ge 0}$ is a $C_0$-semigroup on a Hilbert space $X$, with generator $A$, then $e^{tz} \in \sigma(T(t))\setminus e^{t\sigma(A)}$ implies that $e^{tz}\mathbb T \subseteq \sigma (T(t))$ for almost all $t \ge 0$, where $\mathbb T$ denotes the unit circle. The ``circles conjecture'' is that the result holds more generally for $C_0$-semigroups on Banach spaces, but this remains an open problem. What is known, however, is that if $\T$ is a Banach space $C_0$-semigroup, with generator $A$, and if $\sup \{\|R(a+ib, A)\|: b \in \mathbb R\}=\infty$, then $e^{t a}\mathbb T \subseteq \sigma (T(t))_{t \ge 0}$ for almost all $t \ge 0;$ see \cite[Cor. 4.4]{Wr89}.
If $\R z < s(A)$ for all $z \in \sigma(A)$, then one may set $a=s(A)$ in the inclusion above. It is plausible that the exceptional set here is not merely a null set but necessarily at most countable, as in the Koch-Tataru result.
As an example of an application of Herbst's theorem we mention \cite[Section 4]{BoPe06}, where it was used to clarify the spectral structure of Lax-Phillips semigroups arising in the study of local energy decay.

When $s(A)<\omega(T)$ the identity \eqref{fund} suggests looking more closely at the resolvent of $A$ rather than the spectral bound $s(A)$. This appears to be the right way to proceed on Hilbert spaces and, more generally, on Banach space with good geometry. If $\omega > \omega(T)$ then necessarily $\{z\in\CC: \R z \ge \omega \} \subseteq\rho (A)$ and $\sup_{\R z \ge \omega } \| R(z,A)\| <\infty$. By applying the vector-valued version of Plancherel's theorem on Hilbert spaces and the resolvent identity, one can show that the latter property already characterises the exponential growth bound of $C_0$-semigroups on Hilbert spaces. 

\begin{theorem}\label{gearhart}
If $\T$ is $C_0$-semigroup on a Hilbert space, with generator $A$, then
\begin{equation}\label{omega0}
\omega(T)= \inf \bigg\{\omega > s(A): \sup_{\R z \ge \omega} \|R(z, A)\|<\infty\bigg\} =: s_0 (A).
\end{equation}
 In particular, if $\T$ is bounded, then the following assertions are equivalent:
 \begin{itemize}
  \item[\textup{(i)}] $\T$ is exponentially stable;
  \item[\textup{(ii)}] $i\RR \subseteq\rho (A)$ and $\sup_{s\in\RR} \| R(is,A)\| <\infty$.
 \end{itemize}
\end{theorem}

This result is usually referred to as the Gearhart--Pr\"uss theorem, although it was also obtained independently by Huang \cite{Hg85} and  Monauni \cite{Mo80}; see also \cite{Mo78a}.  It is also a direct consequence of the description of the resolvent set  of $T(t)$ by means of the resolvent of $A$. Indeed, $\mu \in \rho(T(t))\setminus\{0\}$ if and only if the set $\{z\in\CC: e^{zt}=\mu\}\subseteq\rho(A)$ and $R(z,A)$ is bounded on this set;  see e.g \cite{Pr84}. (There is also a Banach space version of the latter result, in which one replaces boundedness of the resolvent by convergence of certain Ces\`aro averages or by an appropriate Fourier multiplier condition; see e.g. \cite{LaRa05}.) Note that the same description of the spectrum of $T(t)$ was obtained, independently and almost simultaneously, by Herbst \cite{He83spec} and Howland \cite{Hw84}. However, these authors did not explore its applications to stability theory.

It is instructive to observe that Theorem~\ref{gearhart} does not generalize to $L^p$-spaces for $p \neq 2$; see for instance \cite[Examples 5.1.11
and 5.2.2]{ABHN01}. Moreover, returning to the Hilbert space setting, the local version of Theorem~\ref{gearhart} does not hold. Indeed, let $\T$ be the (positive) $C_0$-semigroup on $X=L^2(1,\infty)$ given by $(T(t)f)(s)=f(e^t s)$, with generator $A$, and let $f \in X$ be the characteristic function of the set $\bigcup_{n \in \mathbb N}(e^n, e^n +n^{-2})$. Then  the map  $z \mapsto R(z,A)f$ extends to a bounded holomorphic function on $\mathbb C_{-1}$, but the map $t \mapsto \| T(t)f\|$ grows exponentially; see e.g.\  \cite[Example 5.2.3]{ABHN01} for further details. There is a partial remedy for $x \in D((-A)^\alpha)$, where  $\alpha >0$ is fixed. The exponential growth bound for $\|T(t)x\|$, uniform in $x \in D((-A)^\alpha)$, is given by 
\[
\omega_\alpha(T): =\inf \big\{\omega \in \mathbb R: e^{-\omega t}\|T(t)R(\mu, A)^\alpha\|<\infty\big \},
\]
where $\mu \in \rho(A)$ is fixed, and $\omega_\alpha(T)$ can be characterised in terms of  polynomial bounds on $R(z,A)$
 on an appropriate right half-plane, namely
\begin{equation}\label{weiswrobel}
\omega_\alpha(T)= \inf \bigg\{\omega > s(A): \sup_{\R z \ge \omega} (1+|z|)^{-\alpha} \|R(z, A)\|<\infty\bigg\} =:s_\alpha (A) ;
\end{equation}
see \cite{WeWr96}. Moreover, the function $\alpha \mapsto \omega_\alpha(T)$ is convex on $[0,\infty)$. It remains an open question whether it is possible to say more about exponential stability of individual orbits for special classes of semigroups such as $C_0$-semigroups of contractions on Hilbert spaces.

There are important situations in which the spectrum of the generator determines the exponential growth bound of a $C_0$-semigroup, and there exists an intricate abstract theory relating to this question. One such situation arises when the system of root vectors of the generator forms a Riesz basis, which occurs frequently in applications. To our knowledge, the following result by Miloslavskii \cite[Theorem~1]{Ml85} 
is one of the most general results of its kind.

\begin{theorem}\label{Milo}
 Let $(T(t))_{t \ge 0}$ be a $C_0$-semigroup on a Hilbert space $X$, with generator $A$ satisfying
\[
\sigma(A) \subseteq \{z\in\CC : \R z<0\}, \qquad \sigma (A)=\bigcup_{k=1}^{\infty}\Omega_k, 
\qquad \Omega_l \cap \Omega_m=\varnothing, \,\, l \neq m, 
\]
where $\Omega_k$ is a finite set for each $k$. Assume further that the spectral projections $P_k$ corresponding to $\Omega_k$ have finite rank and that the subspaces $X_k:=P_k(X)$, $k \ge 1$, form a Riesz basis of subspaces in $X$. If there exists $\gamma >0$ such that
\[
\sup_{z \in \Omega_k} |\R z| \ge  \gamma  \dim X_k \log (\dim X_k),  
\]
then $(T(t))_{t \ge 0}$ is exponentially stable. If, moreover, $N:=\sup_{k \ge 1} \dim X_k $ is finite, then we have
the explicit estimate
\[
\|T(t)\|\le C (1+t)^{N-1}e^{s(A)t}, \qquad t \ge 0,
\]
for some constant $C>0$.
\end{theorem}

Zabczyk's example mentioned above shows that both parts of this result become false if the relevant assumptions on the dimensions of the spaces $X_k$ are omitted. Note also that  $s(A)=\omega(T)$ whenever $N$ is finite. It is plausible that in the setting of Theorem~\ref{Milo} the WSMT holds for $(T(t))_{t \ge 0}$ if and only if $N$ is finite, but this does not appear to be known.
Note that there exist many  statements in the literature related to Theorem~\ref{Milo}; we mention only the recent papers \cite{AmDiZe14,AmDiZe16}, where the equality $s(A)=\omega(T)$ is proved for generators $A$ whose  eigenvectors form a Riesz basis.  

In many applications, for example those dealing with linearisations of non-linear PDEs, one encounters the problem of finding sharp bounds for the constant $M=M(\omega)$ in the estimate $\|T(t)\|\le Me^{\omega t}$, $t \ge 0$, for a fixed  $\omega\in\RR$, especially when $\T$ is known to be exponentially stable. One thus arrives at the task of obtaining an explicit estimate for $\sup\{e^{-\omega t} \| T(t)\|:t\ge0\}$. It is  usually straightforward to obtain a rough exponential estimate of the form
$
\| T(t)\|\le Le^{\lambda t}$, $t \ge 0.
$
 Helffer and Sj\"ostrand in \cite[Proposition 2.1]{HeSj09} obtained the following useful bound which, incidentally, leads to yet another proof of Theorem~\ref{gearhart}. 

 \begin{theorem}  \label{thm.helffer.sjoestrand}
Let $\T$  be a $C_0$-semigroup on a Hilbert space, with generator $A$. 
Let $L>0$ and $\lambda\in\RR$ be as above and suppose that $\omega < \lambda$
is such that $N(\omega ) :=\sup_{\R z \ge \omega} \|R(z, A)\|$ is finite. Then
\begin{equation}\label{helf}
\|T(t)\|\le L(1 + 2 L N(\omega)\, (\lambda-\omega))e^{\omega t}, \qquad t \ge 0.
\end{equation}
\end{theorem}

It is natural to ask what happens in the limit as $\omega \to \omega(T)$. The answer depends on the rate of blow-up of $R(z, A)$ or, more precisely, of $N(\omega )$ near $\omega (T)$. Assuming that $\omega (T)>-\infty$, $N(\omega)$ blows up as $\omega \to \omega(T)$. If there are $C$, $k>0$ such that $N(\omega) \le C (\omega - \omega(T))^{-k}$ for $\omega - \omega (T)>0$ sufficiently small, then, as shown in \cite{HeSj09}, one gets  $\sup_{t \ge 1} t^{-k} e^{-\omega (T)t}\|T(t)\| <\infty$. 
In fact, one may associate with every blow-up rate of the resolvent an appropriate correction of exponential growth of the semigroup; see \cite{RoVe18a} for details. Recently, using the techniques of \cite{HeSj09}, the following interesting result was proved by Wei in \cite[Theorem~1.3]{We19}. 

\begin{theorem}
Let $\T$ be a $C_0$-semigroup of contractions on a Hilbert space $X$. If $i\mathbb R \subseteq \rho(A)$ and $\omega_0:=\sup_{s\in \mathbb R}\|R(is, A)\|<\infty$, then
\begin{equation}
\|T(t)\|\le e^{\pi/2}e^{-\omega_0 t}, \qquad t \ge 0.
\end{equation}
\end{theorem}

Some Banach space counterparts of Theorem \ref{thm.helffer.sjoestrand} (as well as additional references) may be found in \cite{LaYu13,RoVe18a}, while several interesting and relevant examples of applications are considered e.g.\ in \cite{HeSj09,We19} and \cite[Remark 3.13]{IbMaMa19}.

It is clear that if a $C_0$-semigroup $\T$ is exponentially stable then its orbits $T(\cdot )x$ ($x \in X$) and its weak orbits  $\langle T(\cdot ) x, y\rangle$ ($x\in X$, $y \in X^*$) belong to (vector-valued) $L^p$ and many other classical Banach function spaces on $\mathbb R_+$. A line of theorems originating from Datko's famous result (which are of importance in some applications, e.g.\ in control theory) show that certain converse statements hold, too. We thus have a phenomenon which is Tauberian in character. Moreover, it is possible to estimate the exponential growth bound of a semigroup, as the following result shows. The first part of the result is due to Datko and Pazy (see e.g.\ \cite{Ne96}), and the second part was proved by Weiss \cite{Ws88}.

\begin{theorem}\label{datko}
 Let $\T$ be a $C_0$-semigroup on a Banach space $X$. 
Assume either that there exist $p\in[1, \infty)$ and $C > 0$ such that
\begin{equation}\label{datko_c}
\int_{0}^{\infty} \|T(t)x\|^p\, dt \le C \| x \|^p, \qquad x \in X,
\end{equation}
or, 
 if $X$ is a Hilbert space, that there exist $p\in[1,\infty)$ and $C > 0$ such that 
\begin{equation}\label{weiss_c}
\int_{0}^{\infty} |\langle T(t) x, y\rangle |^p\, dt \le C \| x\|^p \|y \|^p, \qquad x,y \in X.
\end{equation}
Then $\omega(T) \le -(pC)^{-1}$.
\end{theorem}

Note that the existence of constants $C$ such that \eqref{datko_c} and \eqref{weiss_c} hold would already follow from the fact that all orbits (or weak orbits) of $\T$ lie in an $L^p$-space, by a simple application of the closed graph theorem. Observe also that the estimate $\omega(T)\le -(pC)^{-1}$ is optimal in the sense that $\omega(T)$ equals the infimum of the numbers $-(pC)^{-1}$, with $C > 0$ running over all constants for which \eqref{datko_c} holds with respect to an equivalent norm; see \cite[p.\ 82]{Ne96}.

For other results of this type and relevant ideas we refer the reader to \cite[Chapter~3]{Ne96}, \cite{ABHN01}, and to \cite{Ni02}. There are many more statements along the lines of Theorem~\ref{datko} in the literature (although most of them are primarily of theoretical interest).

While there are many perturbation results for the exponential stability, most of them concern perturbations which are small in a metric sense; in this case exponential stability of the perturbed semigroup comes as no surprise. However, there are some nice exceptions of this rule, as the following recent result by Pr\"uss \cite{Pr15}  shows.

\begin{theorem}
If $A$ generates an exponentially stable $C_0$-semigroup $(T(t))_{ t \ge 0}$ on a Banach space $X$, and if $B \in {\mathcal L} (X)$ is such that $\lim_{t\to 0+}\|(T(t)-I)B\|=0$ (for example, if $B$ is compact) and $\lambda - A - B$ is invertible for $\lambda \in \overline{\mathbb C_+}$, then $A+B$, too, generates an exponentially stable $C_0$-semigroup.
\end{theorem}


\subsection{Strong stability}
Strong stability is another basic stability concept in the theory of operator semigroups. Unlike exponential stability, it is distinctly qualitative in character.
Nevertheless, it is precisely this notion of stability that lies behind many of the more recent developments in the quantified asymptotic behaviour of operator semigroups.

Recall that a $C_0$-semigroup $\T$ on a Banach space $X$ is said to be {\em strongly stable}, or simply {\em stable}, if for every $x\in X$,
\[
 \lim_{t\to\infty} \| T(t)x\| = 0.
\]
By the uniform boundedness principle, every stable semigroup is necessarily bounded and, by the spectral inclusion theorem, the generator $A$ of a stable semigroup satisfies $\sigma (A) \subseteq \{z\in\CC: \R z \leq 0 \}$. The example of the left-shift semigroup on $L^2 (\RR_+ )$ shows that a semigroup may be stable and yet have a generator $A$ whose spectrum is ``maximal'' in the sense that $\sigma (A) = \{ z\in\CC:\R z \leq 0\}$. Thus the location of the spectrum of the generator alone does not tell us a great deal about stability of the semigroup. 
The initial intuition behind the spectral approach is  that the smaller the spectrum of the generator on the imaginary axis the better the stability properties of the associated semigroup. A particularly revealing illustration of this is the following famous result due to Arendt and Batty, and Lyubich and V\~u \cite[Theorem 5.5.5]{ABHN01}.

\begin{theorem}[Arendt-Batty-Lyubich-V\~u]\label{arendt}
Let $\T$ be a bounded $C_0$-semigroup on a Banach space $A$, and suppose that the generator $A$ is such that the boundary spectrum $\sigma(A)\cap i\RR$ is at most countable and contains no eigenvalues of the adjoint of $A$. Then $\T$ is stable.
\end{theorem}

The result is best possible as far as the spectral assumptions on $A$ are concerned; see e.g. \cite{ABHN01}. Subsequently, more sophisticated results of this kind, applicable in the case when the spectrum contains arbitrary large sets on the imaginary axis, were proved in the series of papers  \cite{To01,ChTo03,ChTo04, BoChTo07}; see also the survey \cite{ChTo07}. Here the above ``spectral smallness'' principle was replaced by a more general one: the more the \emph{resolvent} of the generator can be extended to or across the imaginary axis, the better the stability properties of the semigroup orbits. To convey the flavour of the relevant results we formulate a statement which turned out to be particularly useful in some applications to PDEs; see \cite[Theorem~5]{BoLe95}, \cite[p. 75-76]{To01} and \cite{ChTo07} for further details.

\begin{theorem}
Let $\T$ be a completely non-unitary $C_0$-semigroup of contractions on a Hilbert space $X$,
with generator $A$.
Then $\T$ is stable if and only if the set
\[
\Big\{x \in X: \lim_{a \to 0+} \sqrt{a}R(a+ib, A)x=0  \text{ for almost all } b\in\RR\Big\}
\]
is dense in $X$. If $\T$ is merely assumed to be bounded then it is stable provided the set
\[
\Big\{x \in X: \lim_{a \to 0+} \sqrt{a} R(a+ib, A)x=0  \text{ for all } b\in\RR\Big\}
\]
is dense in $X$.
\end{theorem}

Note that it is possible to formulate necessary and sufficient (integral) conditions for stability of a bounded $C_0$-semigroup $\T$ or its individual orbits. For an account of these results concerning strong stability, including stability of semigroups on Banach spaces,  characterisations of stability in non-resolvent terms, and even stability of (non-autonomous) evolution families,   we refer the reader to the survey \cite{ChTo07} and the references therein.

One may also define a notion of weak stability for $C_0$-semigroups. Indeed, one says that $(T(t))_{t \ge 0}$ is \emph{weakly stable} if all of the orbits of $(T(t))_{t \ge 0}$ converge to zero in the weak topology as $t \to \infty$. The notion of weak stability corresponds to the concept of mixing in ergodic theory, where one deals with weakly stable unitary $C_0$-groups of Koopman operators on $L^2$-spaces. The general theory of weak stability of $C_0$-semigroups is not particularly well developed, and in particular weakly mixing unitary groups
are still studied in ergodic theory largely on a case-by-case basis.
On the other hand, weak stability can often serve as a Tauberian condition for strong stability, as the following statement shows.

\begin{theorem}\label{weak_stab}
Let $\T$ be a $C_0$-semigroup of contractions on a Hilbert space $X$, with generator $A$. 
If $\T$ is completely non-unitary then it is weakly stable. Moreover,
if $A$ has compact resolvent then the following conditions are equivalent:
\begin{itemize}
\item [(i)] $\T$ is weakly stable;
\item [(ii)] $\T$ is strongly stable;
\item  [(iii)] $(T^*(t)_{t\ge0}$ is strongly stable;
\item [(iv)] $\T$ is completely non-unitary.
\end{itemize}
If $A=A_0-BB^*$, where $A_0$ is skew-adjoint, has compact resolvent, and $B$ is bounded, then each of the conditions (i)-(iv)
is equivalent to the property that $B^* x\neq0$ for every non-zero $x \in X$ such that $A_0x=is x$ for some $s \in \mathbb R$.
\end{theorem}

This result has appeared, in one form or another, in a number of papers. It is used implicitly in the paper \cite{Iw69}, which studies energy decay for a class of hyperbolic equations including the damped wave equation; see \cite{BaLeRa89} for a discussion. Proofs of the equivalence of (i)-(iv)  in Theorem~\ref{weak_stab} may be found in \cite{Bc78,Lv80,LuGuMo99}, while the last statement of the theorem is contained in \cite[Theorem~5]{KoSk84} and \cite[Theorem~14]{BaVu90}. As another statement  similar in flavour to Theorem~\ref{weak_stab}, recall that weak stability of a bounded semigroup $\T$ on an $L^1$-space implies strong stability (see \cite{ChTo07}), and the same holds on $C(K)$-spaces, where $K$ is a compact Hausdorff space, if $\T$, in addition, is irreducible; see \cite{MeTr93}. For further results in which the stability (or, more generally, convergence) of a $C_0$-semigroup is a consequence of mild assumptions on the geometry of the underlying Banach space, see \cite{Gl16b} and the references therein. 

\section{Quantified Tauberian theorems and semi-uniform stability of operator semigroups}\label{sec:semi}

In many cases, stability results for $C_0$-semigroups or their individual orbits are consequences of general Tauberian theorems for Laplace transforms. The latter often serve as an inspiration for the former, and in modern treatments both are recognised as being two sides of the same coin. Nowadays, the role of Tauberian theory in the theory of $C_0$-semigroups is well understood. The relevant theory was developed in a number of papers, arguably starting with the seminal paper \cite{ArBa88}. Later, in \cite{Ch98}, a finer distributional approach was developed,  which has recently been extended to include quantitative aspects; see  \cite{ChSe16, BaBoTo16, DeSe19,DeSe19a}. A good introduction to modern Tauberian theory may be found in \cite{Ko04}, while applications to operator semigroups are discussed thoroughly in  \cite{ABHN01}.

The following result is usually attributed to Ingham \cite{In35} and Karamata \cite{Ka34}, although a version was in fact first discovered by M. Riesz \cite[Satz II]{Ri1924}. It is a classical Tauberian theorem which underpins just about all of what follows here; we refer to \cite{ChTo07} for a short proof.

\begin{theorem}\label{ingham}
Let $X$ be a Banach space and let $f :\mathbb R_+ \to X$ be a bounded and uniformly continuous function whose Laplace transform $\widehat f$ admits a holomorphic extension to each point of $i\mathbb R$.
Then $\lim_{t\to\infty} \| f(t)\| = 0$.
\end{theorem}

Theorem~\ref{ingham} is the starting point for numerous results in the asymptotic theory of $C_0$-semigroups, and it has been generalised in various directions. It is natural in analysis, and also in Tauberian theory, to equip a convergence result with a rate whenever this is possible. In the case of Theorem~\ref{ingham} this can be done in a way which has a number of interesting consequences, not least for operator semigroups. Observe, however, that Theorem~\ref{ingham} is best possible in the sense that one cannot expect any rate of decay for $f$ if no further assumptions are imposed on the growth of $\widehat f$, even if $\widehat f$ extends to an entire function; see e.g.\ \cite{DeVi18,DeSe19,BrDeVi20}. In order to obtain quantitative results one assumes that $\widehat f$ extends analytically beyond $i\mathbb R$ to some precise domain and that this analytic extension satisfies an appropriate bound in this domain.
Given  continuous and increasing functions $M$, $K: \mathbb R_+\to (0,\infty)$ we define
\begin{equation}\label{mlog}
\Omega_M:=\left \{z \in \mathbb C: \R z > - \frac{1}{M(|\I z|)}\right \} \text{ and } 
M_{K}(s) := M(s) \big(\log(1 +K(s)) + \log (1 + s)\big)
\end{equation}
for $s\ge0$. Note that the function  $M_{K}$ is continuous, strictly increasing and unbounded, and hence it has an inverse function $M_{K}^{-1}$ defined on $[a,\infty)$ for some $a>0$. The function $M$ itself may not be strictly increasing, and we denote by $M^{-1}$ any right-inverse of $M$. Without further comment we extend both inverses by $0$ to $\RR_+$. The following result is a rather general Tauberian theorem. 

\begin{theorem}  \label{x0}
Let $f \in L^1_{\loc}(\RR_+;X)$ be a function with weak derivative $f'\in L^p (\RR_+ ;X)$, where $p\in [1,\infty ]$. Let $M$, $K: \mathbb R_+\to (0,\infty)$ be continuous and increasing functions satisfying, for some $\varepsilon \in (0,1)$, $C\geq 0$, 
\[
 M(s) \leq K(s) \leq C\, e^{e^{(sM(s))^{1-\varepsilon}}}, \qquad s\in\RR_+ . 
\]
Assume that $\widehat f$ extends analytically to $\Omega_M$ and that
\begin{equation} \label{resboundp}
\|\widehat f (z)\| \le K( |\I z|) , \qquad z\in \Omega_M .
\end{equation}
Then there exists $c>0$ such that 
\begin{equation} \label{Mlogestp}
( t \mapsto M_{K}^{-1} (ct) \| f(t)\| ) \in L^p (\RR_+ )  . 
\end{equation}
\end{theorem}

This Tauberian theorem is the result of developments over a period of almost ten years. In the special case $p=\infty$ and $K=M$ it was first obtained by Batty and Duyckaerts in \cite{BaDu08}. In that case, the function $M_K (s) = M(s) ( \log (1+M(s)) + \log (1+s))$  is often denoted by $\Mlog$. It is the achievement of Batty, Borichev and Tomilov \cite{BaBoTo16} to extend the basic result from \cite{BaDu08} and the classical Ingham-Karamata-Riesz theorem to general $p$ and to allow functions $K$ of the form $(1+s)^\alpha M(s)^\beta$ (for some $\alpha$, $\beta\geq 0$) instead of just $K=M$. In this case $M_K$ and  $M_{\log}$ do not differ essentially. There is also an extensive discussion in \cite[Section 4]{BaBoTo16} of the fact that more general functions $K$ are possible, too. In fact, there are sometimes non-trivial relations between the functions $M$ and $K$, 
and the interplay between the shape of $\Omega_M$ and the growth of $\widehat f$ in $\Omega_M$ can be analysed by means of more or less standard techniques of complex analysis, such as the theory of harmonic measure. For instance, if $M(s)=C(1+s)^\alpha$  and if $\widehat f$ grows polynomially in $\Omega_{M}$, then, as was noted  in \cite[Lemma 3.4]{BoTo10}, one necessarily has $|\widehat f(z)|\le C_\epsilon(1+|z|^{\alpha +\epsilon})$ for all $z\in\Omega_{cM}$ for an appropriate constant $c>1$. Theorem \ref{x0} in its full generality was finally proved in an unpublished manuscript by Stahn \cite[Theorem 1.1]{Sta17a}; see also \cite{Sta18a}. 

Regardless of $p$, Theorem~\ref{x0} is optimal in the sense that one cannot obtain any better weights in \eqref{Mlogestp} than $M_{K}^{-1}$. Optimality for $p=\infty$ and polynomial $M = K$ was proved in \cite[Theorem~3.8]{BoTo10} by an explicit and fairly involved construction.  By developing further the ideas from \cite{BoTo10}, this approach was then extended to all $p$ and to some other functions $M$ (and also to individual orbits of Hilbert space semigroups) in \cite[Sections 5 and 7]{BaBoTo16}, and in \cite[Theorem 4.1]{Sta17a} for an even larger class of functions $M$ and $K$; see also \cite{DeSe19} and \cite{DeSe19a} for an alternative abstract approach to optimality.

In order to make the link between general Tauberian theory and $C_0$-semigroups,
it suffices to observe that if $(T(t))_{t \ge 0}$ is a bounded $C_0$-semigroup on a Banach space $X$ with generator $A$, and if $A$ has no spectrum on the imaginary axis, then applying Theorem~\ref{ingham} to individual orbits shows that the semigroup is stable. This is of course a very special case of Theorem~\ref{arendt}. Applying Theorem~\ref{ingham} to the bounded, Lipschitz continuous, operator-valued function $T(\cdot )A^{-1}$ yields  
\begin{equation}\label{semi}
\lim_{t \to \infty} \|T(t)A^{-1}\|=0.
\end{equation}
To distinguish this kind of asymptotic behaviour of $\T$, we call a bounded $C_0$-semigroup $\T$ satisfying \eqref{semi}  {\em semi-uniformly stable}, noting  that in this case the stability of $\T$ is uniform with respect to initial values from the unit ball of $D(A)$ endowed with the graph norm. We emphasise that our terminology is not particularly common in the literature. However, the terminology is at least natural from the point of view that every exponentially stable semigroup is semi-uniformly stable, and every semi-uniformly stable semigroup is (strongly) stable. It is straightforward to see that the converses of these statements are both false.

A natural question, both from an operator-theoretic perspective and from the point of view of applications to PDEs, is whether semi-uniform stability can be quantified.
To clarify the quantification problem note first that convergence of $\|T(\cdot)A^{-1}\|$ to zero  can formally  be translated into an estimate of the form
\begin{equation} \label{eq.decay.rate}
\| T(t)x\| = O( r(t)),\quad t\to\infty,
\end{equation} 
for some function $r\in C_0 (\RR_+)$ and for every $x\in D(A)$, that is, for every {\em classical} solution of the Cauchy problem $\dot u = Au$. If the semigroup $\T$ is stable but not exponentially stable, then one can never expect such an estimate to hold for {\em every} initial value $x\in X$. Indeed, having \eqref{eq.decay.rate} for every $x\in X$ is equivalent  to  exponential stability by a simple application of the uniform boundedness principle, as has already been observed. In fact, for strongly stable semigroups which are not exponentially stable we obtain the following statement; see \cite{Ne96,MuTo13}.

\begin{theorem}\label{lower_b}
\begin{itemize}
\item [(a)] Let $\T$ be a strongly stable but not exponentially stable $C_0$-semigroup on a Banach space $X$. Then for every $r\in C_0 (\RR_+)$ there exists $x\in X$ such that $\| T(t) x\| \geq r(t)$ for all $t\geq 0$. 
\item [(b)] Let $\T$ be a weakly stable but not exponentially stable $C_0$-semigroup on a Hilbert space $X$. Then for every $r\in C_0 (\RR_+)$ there exists $x\in X$ such that $|\langle T(t)x,x\rangle| \ge r(t)$ for all $t \ge 0$.
\end{itemize}
\end{theorem}

Note that for non-exponentially stable semigroups the above statements also rule out the possibility of various integral conditions being satisfied for all semigroup orbits. This is in accordance with Theorem~\ref{datko} stated above.
For more results on lower bounds for weak orbits of $C_0$-semigroups on mostly reflexive Banach spaces, see \cite{MuTo13} and \cite{Sto10}. 

While the decay of a stable semigroup can be abitrarily slow, the decay of sufficiently regular orbits admits a rate whenever the semigroup is semi-uniformly stable. This claim is made precise by the following theorem, a substantial part of which can be deduced from the Tauberian theorem \ref{x0}. The statement was obtained in \cite{BaDu08} by Batty and Duyckaerts, as a culmination of several notable preceding results originating from both abstract operator theory and its applications, e.g.\ in \cite{Le96,Bq98,BEPS06,LiRa05}. The paper \cite{BaDu08} sparked a considerable amount of further research on quantitative aspects of semi-uniform stability for $C_0$-semigroups, and it is a precursor of many of the quantitative stability results discussed below.

\begin{theorem}\label{semigrouprates} 
Let$\T$ be a bounded $C_0$-semigroup on a Banach space $X$, with generator $A$, and let $\mu \in \rho(A)$. Then the following are equivalent:
\begin{enumerate}
\item [{(i)}] $\sigma(A)\cap i\mathbb R$ is empty; \label{sgr1}
\item [{(ii)}] $\lim_{t \to \infty} \|T(t)R(\mu,A)\|=0$.
 \end{enumerate}
Moreover, if either (i) or (ii) hold, let
\begin{equation*}
M(s) := \sup \{\|R(ir,A)\| : |r| \le s\}, \qquad s\ge 0,
\end{equation*}
define the function $M_{\log}$ by $M_{\log} = M_K$ where $K=M$ and $M_K$ is as in \eqref{mlog}, and  let $M^{-1}$ denote any right-inverse of $M$. Then there exist positive constants $C$, $C'$, $c$, $c'$ such that
\begin{equation} \label{genestintro}
\frac{C'} {M^{-1}(c't)} \le \|T(t)A^{-1}\| \le \frac{C}
{M_{\log}^{-1}(ct)}
\end{equation}
for all sufficiently large $t$.
\end{theorem}

The implication (i)$\implies$(ii) is proved for instance in \cite[p. 803]{Ba90} and \cite[Corollary~3.3]{Vu92}, although it was already contained implicitly in \cite{ArBa88}; it may be viewed as a consequence of the Ingham-Karamata-Riesz theorem. In order to prove the converse implication, assume that (ii) holds and that $s\in \RR$ is such that $is\in\sigma(A)$. Then $e^{ist}(1-is)\inv\in \sigma(T(t)R(1,A))$  by the spectral inclusion theorem \cite[Theorem~16.3.5]{HiPh57} for the Hille-Phillips functional calculus, and in particular $|1-is|\inv\le \|T(t)R(1,A)\|$ for all $t\ge 0$. Letting $t\to\infty$  yields the desired contradiction, since  (ii) in particular holds for $\mu=1$ by a simple application of the resolvent identity.

Observe that if $\|R(is, A)\| \le M(|s|)$ for $s \in \mathbb R$, then the Neumann series expansion shows that, for any $c > 1$, the region $\Omega_{cM}$ is contained in the resolvent set of $A$  and $\|R(z, A)\|\le \frac{c}{c-1} M(|\I z|)$ for $z \in \Omega_{cM}$.  Hence the upper bound in \eqref{genestintro} is a direct consequence of Theorem~\ref{x0}. In contrast to the setting of Theorem~\ref{x0}, however, in this case it is natural to describe the region of holomorphic  extension of $R(z, A)$ and the growth of the extension by the same function $M$, up to constant multiples. The lower bound in \eqref{genestintro} arises as consequence of simple operator-theoretic considerations based on the semigroup analogue of the fundamental theorem of calculus.

Note that if $M(s)=Ce^{as}$, then $M^{-1}$ and $M_{\log}^{-1}$ are asymptotically of the same order, and Theorem~\ref{semigrouprates} yields the optimal logarithmic rate of decay for $\|T(t)A^{-1}\|$. This result, which goes back  to \cite{Le96} and, with an optimal rate, to Burq \cite{Bq98}, is of great value for applications, for example to the study of damped wave equations. By going to the opposite end of the scale and taking a constant function $M$ one recovers part of the Weis-Wrobel relation \eqref{weiswrobel} for bounded $C_0$-semigroups in the case  $\alpha=1$.  A general local version of this statement is due to van Neerven \cite{Ne96}; see also \cite{Ba03} the references therein. Thus it is natural and important when the upper bound in \eqref{genestintro} matches the lower bound. It can be easily checked that for semigroups of normal operators and for sufficiently regular $M$ one can indeed replace the upper bound $1/M_{\log}^{-1}$  by the lower bound $1/M^{-1}$; see e.g.\ \cite[Section 5.1]{BaChTo16}. However, it is known that the function $M_{\log}^{-1}$ in \eqref{genestintro} cannot in general be replaced by $M^{-1}$ for bounded \emph{Banach space} $C_0$-semigroups. As in the case of Theorem~\ref{x0}, Theorem~\ref{semigrouprates} is optimal in the sense that one cannot obtain better upper estimates for the semigroup orbits than the one given by \eqref{genestintro}. In the case of polynomial resolvent growth this goes back once again to \cite{BoTo10}, but more recently an elegant abstract approach for obtaining optimal lower bounds for semigroup orbits was found by Debruyne and Seifert in \cite{DeSe19,DeSe19a}. Note, however, that the lower bounds contained in these optimality results hold only along a subsequence of $\mathbb R_+$. It is not known whether one may prove stronger optimality results in the spirit of Theorem~\ref{lower_b}.  

Concerning the closely related situation of so-called cut-off (or Lax-Phillips) semigroups it is important to note that around the same time as Batty and Duyckaerts in \cite{BaDu08}, Christianson proved a very similar decay rate result for functions of the form $f (t)= \chi_1 T(t) R(1,A)^{k}\chi_2$, where $\T$ is a $C_0$-semigroup of contractions on a Hilbert space, $\chi_1$ and $\chi_2$ are bounded linear operators, the Laplace transform $\widehat{f}$ extends to a domain $\Omega_M$ where it satisfies a polynomial growth estimate with $K(s)=C(1+s)^N$, and $k\geq N+2$ \cite[Theorem~3]{Cr09}. The result states that $\| f(t)\| = O(1/\tilde{M}_{\log}^{-1} (ct)^{k/2})$, where $\tilde{M}_{\log} (s) = M(s) \log (1+s)$. This result has been applied to the study of concrete equations, in particular damped wave equations, for instance in \cite{Cr09,CSVW14,Ri14,Jin20}. It is worth pointing out that, unlike Batty and Duyckaerts, Christianson in his assumptions already uses separate functions $M$ and (a polynomial) $K$ for the description of the region of analytic continuation and for the estimate of the cut-off resolvent $ \chi_1 R(z ,A)R(1,A)^k \chi_2$ in that region. However, improved versions of the Tauberian Theorem \ref{x0} (see \cite[Theorem 4.2]{BaDu08}, \cite[Theorem 1.1]{Sta17a}) yield the decay rate $\| f(t)\| = O(1/M_{\log}^{-1} (ct)^k)$; see also \cite[Corollary 4.2]{BaDu08} which is especially designed for the situation of cut-off semigroups. Note that $M_{\log}^{-1}$ and $\tilde{M}_{\log}^{-1}$ are comparable for polynomial $M$, and that therefore the result from \cite{BaDu08} is stronger in this case.

Of course, from the point of view of applications to PDEs, one's primary interest is in the situation where $\T$ is a bounded $C_0$-semigroup on a Hilbert space. 
It turns out that for such semigroups Theorem~\ref{semigrouprates} can be substantially sharpened,  by replacing $M_{\log}$ with $M$  for a very large class of functions $M$, thus answering an open question in \cite{BaDu08} for polynomial $M$. As in the situation of the Gearhart-Pr\" uss theorem, the availability of the vector-valued Plancherel theorem plays a crucial role here. The following result was first obtained in \cite[Theorem~2.4]{BoTo10}. It shows that if, in Theorem~\ref{x0}, $X$ is a Hilbert space and $M$ is of polynomial form, then the upper bound in \eqref{genestintro} can be sharpened to match the lower bound, and in particular the sharpened  upper bound is optimal. 

\begin{theorem}\label{thm:BT} 
Let $\T$ be a bounded $C_0$-semigroup $\T$ on a Hilbert space, with generator $A$, and let $\alpha>0$ be a constant. The following assertions are equivalent:
\begin{itemize}
 \item[\textup{(i)}] $\|T(t)(I-A)^{-1}\|=O(t^{-1/\alpha})$ as $t\to\infty$;
 \item[\textup{(ii)}] $i\RR\subseteq\rho (A)$ and $\|R(is,A)\| = O(|s|^\alpha)$.
\end{itemize}
\end{theorem}

Note that  polynomial decay rates as in Theorem~\ref{thm:BT} are stable with respect to certain classes of factored perturbations $BC$ of $A$, where the operators $B$ and $C$ are subordinated to appropriate fractional powers of $-A$ and $-A^*; $ see e.g.\ \cite{Pa12} for details.

The approach of \cite{BoTo10} was refined and extended in \cite{Ma11,BaChTo16}. 
Here the study of uniform stability of $T(\cdot)x$ for $x \in D(A)$ was extended to  orbits  with $x \in \I (A)$ and also,  as a combination of the two,
to uniform stability of orbits with  $x \in \I (A)\cap D(A)$. Uniform stability of $T(\cdot)x$ for $x \in \I (A)$ corresponds in operator-theoretic terms to $\|T(t)A(I-A)^{-1}\|\to 0$, and the latter condition is characterised (in the Hilbert space setting) by the spectral condition $\sigma(A)\cap i\RR\subseteq \{0\}$ together with the resolvent condition $\sup_{|s|\ge 1}\|R(is,A)\|<\infty$. 
The spectral condition on its own characterises (even on Banach spaces) the decay $\|T(t)A(I-A)^{-2}\|\to 0$, which is the operator-theoretic formulation of uniform stability of orbits $T(\cdot)x$ for $x\in \I (A)\cap D(A)$. 
Apart from providing a general framework within which to analyse these three types of semi-uniform stability, the paper \cite{BaChTo16} also presented a number of quantified (and often optimal) results on semi-uniform stability, which are in the the spirit of Theorem~\ref{thm:BT} but allow for resolvent growth which is not necessarily polynomial but more generally regularly varying. This direction of research was further developed for instance in in \cite{Sf15,ChSe16,RoSeSt19}, where a number of generalisations and complementary results may be found.

The most notable extension of the results in \cite{BoTo10,BaChTo16} was obtained in \cite{RoSeSt19} by Rozendaal, Seifert and Stahn.  A continuous increasing function $M: \RR_+\to(0,\infty)$ is said to be  of \emph{positive increase} if there exist positive constants  $\alpha>0$, $c\in(0,1]$ and $s_0> 0$ such that 
\[
\frac{M(\lambda s)}{M(s)}\ge c\lambda^\alpha,\quad \lambda\ge1, \,s\ge s_0.
\]
This is a large class of functions, which in particular contains the class of \emph{regularly varying} functions (of positive index) considered in \cite{BaChTo16}.  The following result was proved in \cite[Theorem~3.2]{RoSeSt19}.

\begin{theorem}[Rozendaal-Seifert-Stahn] \label{thm:RSS}
 Let $\T$ be a bounded $C_0$-semigroup on a Hilbert space $X$, with generator $A$, and let $M:\RR_+\to (0,\infty )$ be a function of positive increase. The following assertions are equivalent:
\begin{itemize}
 \item[(i)] $\|T(t)(1-A)^{-1}\|=O(1/M^{-1} (t) )$ as $t\to\infty$;
 \item[(ii)] $i\RR\subseteq\rho (A)$ and $\|R(is,A)\| = O(M(|s|))$ as $|s|\to \infty$.
\end{itemize}
\end{theorem}

Note that, as the analysis of multiplication semigroups in \cite[Section 5.1]{BaChTo16} reveals,  the class of functions of  positive increase is the largest possible class for which the resolvent estimate $\|R(is,A)\|\le M(|s|)$, $s\in\RR$, implies the decay rate $\|T(t)A\inv\|=O(1/M\inv(t))$ as $t\to\infty$; see also \cite{RoSeSt19}. 

It sometimes happens in applications that a resolvent estimate on $i\mathbb R$ is known only along a subsequence. To clarify this situation and complement Theorem~\ref{semigrouprates}, observe that if $M: \mathbb R_+\to(0,\infty)$ is a continuous increasing function such that $M(s)\to\infty$ as $s\to\infty$ and 
$\limsup_{|s|\to\infty} M(|s|)\inv\|R(is,A)\|>0$,
then there exists $c>0$ such that 
$\limsup_{t\to\infty}M\inv(ct)\|T(t)(I-A)^{-1}\|>0$,
and if $M$ has positive increase then the latter inequality holds for \emph{all} $c>0$; see \cite[Proposition 5.4]{CPSST19}.
In fact, the same applies more generally to bounded $C_0$-semigroups on \emph{Banach} spaces.

Exponential and semi-uniform stability with polynomial rates can also be approached by means of Fourier multiplier techniques. We are unable to discuss this approach in our short survey, and we instead refer to the interesting recent paper \cite{RoVe18} (and the references therein) for details. In particular, the paper contains versions of Theorem~\label{thm:BT} for unbounded semigroups; see e. g.\ \cite[Section~4]{RoVe18}.

\section{Applications to second-order Cauchy problems}\label{damped}

Since this survey emphasises abstract aspects we found it instructive, in passing to descriptions of illustrative applications of our techniques,  to begin with a discussion of an abstract model  and then to illustrate briefly, in the concluding remarks, how it relates to several concrete situations dealing mostly with damped wave equations. This subsection is based on the presentation in a very nice paper by Anantharaman  and L\'eautaud \cite[Part II]{AnLe14}, which in fact considers a slightly more general setup than ours. Note that similar arguments can be found in a number of papers and books dealing with asymptotic behaviour of hyperbolic equations and operator matrices, but we refrain from mentioning all of them here. 

For the remainder of this section, let $A$ and $B$ be self-adjoint positive semi-definite operators on a Hilbert space $X$. Assume that $A$ has compact resolvent, and that $B$ is bounded. Assume further that 
\begin{equation} \label{calderon}
 \| Bx\| >0 \text{ for every (non-zero) eigenvector $x$ of } A .  
\end{equation}
The general approach given below in fact works under much milder assumptions on $A$ and $B$, but in view of our applications to damped wave equations it is reasonable to restrict our attention to this particular setting. Consider then the following second-order Cauchy problem in $X$:
\begin{eqnarray}\label{equation}
\ddot u + B \dot u + A u=0, \qquad u(0) = u_0 , \,\, \dot u (0) = u_1 , 
\end{eqnarray}
for certain initial values $u_0$, $u_1\in X$. Setting $U=(u, \dot u)$, this second-order Cauchy problem can be rewritten on the space $\mathcal X := D(A^{1/2}) \times X$ as a first-order Cauchy problem
\begin{eqnarray}\label{Cauchy_f_ord}
\dot U =\mathcal A U, \qquad    
U(0)= (u_0, u_1) \in \mathcal X ,
\end{eqnarray}
where 
\begin{equation}
\mathcal A = \left( \begin{array}{cc} 0 & I \\ -A & -B \end{array} \right), \qquad \text{with} \quad D(\mathcal A)=D(A)\times D(A^{1/2}).
\end{equation}
The operator $\mathcal A$ generates a $C_0$-semigroup of contractions $(\mathcal T(t))_{t \ge 0}$  on $\mathcal X$ by the Lumer-Phillips theorem and a simple perturbation argument. Hence for every pair of initial values $(u_0,u_1) \in D(A)\times D(A^{1/2})$  the second order problem \eqref{equation} admits a unique classical solution 
$
u \in C(\mathbb R_+;
D(A))\cap C^1(\mathbb R_+; D(A^{1/2})) \cap C^2 (\RR_+ ;X)
$.

Since $A$ has compact resolvent, the embeddings $D(A) \hookrightarrow D(A^{1/2})\hookrightarrow X$ 
are compact, and boundedness of $B$ implies that $\mathcal A$ too has compact resolvent. The spectrum of $\mathcal A$ contains therefore only isolated eigenvalues of finite multiplicity. It was proved in \cite{AnLe14} that the spectrum $\sigma (\mathcal A)$ is contained in the union of the open strip $\{z\in\CC: -\frac 12 \|B\| < \R z < 0\}$ and the interval $[-\|B\|,0]$. Moreover, ${\Ker}\, \mathcal A={\Ker}\, A \times \{ 0\}$, so that $0\in\sigma (\mathcal A)$ if and only if $0\in\sigma (A)$. Note that assumption \eqref{calderon} plays a crucial role in showing that $\sigma(\mathcal A)\cap i\RR\subseteq\{0\}$.

If $0\in\sigma (\mathcal A)$ we let $\mathcal P_0$ be the Riesz projection corresponding to the eigenvalue $0$ (noting that $\mathcal P_0$ need \emph{not} be orthogonal), and if $\mathcal A$ is invertible we set $\mathcal P_0 = 0$. Let $\mathcal X_0 := (I-\mathcal P_0) \mathcal X$ and $\mathcal A_0:=\mathcal A|_{\mathcal X_0}$. Then $\mathcal A_0$ is the generator of the $C_0$-semigroup $(\mathcal{T}_0(t))_{t\ge0}$ of contractions (namely the restriction of the semigroup $(\mathcal{T}(t))_{t\ge0}$ to the $\mathcal T$-invariant subspace $\mathcal X_0$) and $\sigma (\mathcal A_0)=\sigma (\mathcal A)\setminus \{0\}$. In particular, it follows from Theorem~\ref{semigrouprates} that the semigroup $(\mathcal{T}_0(t))_{t\ge0}$ is semi-uniformly stable. 

One would typically like to determine whether $(\mathcal{T}_0(t))_{t\ge0}$ is exponentially stable or, if not, to obtain uniform decay rates for classical solutions. The decay of $(\mathcal{T}_0(t))_{t\ge0}$ is intimately related to the decay of the physically more relevant {\em energy} of solutions $u$ of the second-order problem  \eqref{equation}:
\begin{equation}\label{abstract_energy}
E(u,t):=\frac12 \left(\|A^{1/2}u\|^2+\|u_t\|^2\right) .
\end{equation}
In fact, the decay rate of the energy $E(u,t)$ for all pairs of initial values $(u_0,u_1)\in D(A)\times D(A^{1/2})$ (no restriction to $\mathcal X_0$!) is the same as the square of the decay rate of  $\|\mathcal T_0(t)\mathcal A_0^{-1}\|$ and, as a result, we may apply the results in Section \ref{sec:semi}. It is then natural to say that the energy of \eqref{equation} decays semi-uniformly with rate $r\in C_0 (\RR_+ )$ if $E (u,t) \leq C  r(t)^2 \|\mathcal A (u_0, u_1)\|^2$ for all $(u_0,u_1)\in D(A) \times D(A^{1/2})$. Note that $E(u,t) = O( r(t)^2)$ for some $r\in C_0 (\RR_+ )$ and for {\em all} initial values $(u_0,u_1)\in\mathcal X$ is equivalent to exponential stability of the semigroup $(\mathcal{T}_0(t))_{t\ge0})$, and therefore to exponential decay of $E(u,t)$ for {\em all} solutions of \eqref{equation}.

Recall from Section \ref{sec:exp_str} the discussion of  circles contained in the spectrum of $C_0$-semigroups. This discussion implies in particular that the well-known alternative --
either $\|\mathcal T_0 (t)\|=1$ for all $t\ge 0$ or $(\mathcal T_0 (t))_{t \ge 0}$ is exponentially stable -- holds here in a much stronger form: either $(\mathcal T_0 (t))_{t \ge 0}$ is exponentially stable or $\mathbb T \subseteq \sigma (\mathcal T_0 (t))$ for almost every $t \ge 0$. Thus, unless $(\mathcal T_0(t))_{t \ge 0}$ is exponentially stable, its peripheral spectrum (in fact, even the peripheral essential spectrum) is as large as possible for almost all values of $t$.

In order to link the study of $(\mathcal T_0 (t))_{t \ge 0}$ with the results from the preceding sections, define the quadratic operator pencil $P(z):=z^2 + z B + A$, $z \in \mathbb C$, with $D(P(z))=D(A)$. The next statement in combination with the abstract decay criteria from Section \ref{sec:semi} tells us that semi-uniform decay rates for the energy $E$ can be rewritten as operator norm estimates for $P(z)^{-1}$ with $z\in i\RR$; see \cite[Lemma~4.6]{AnLe14}.

\begin{theorem}\label{resolvent_estimate}
There exists $C>1$ such that for all $s\in \mathbb R$ with $|s| \ge 1$, 
\begin{align*}
\|R(is, \mathcal A_0)\|-C|s|^{-1}\le& \|R(is, \mathcal A)\|\le \| R(is, \mathcal A_0)\|+ C |s|^{-1}, \text{ and}\\
C^{-1}|s|\|P(is)^{-1}\| \le& \|R(is, \mathcal A)\|\le C(1+|s|)\|P(is)^{-1}\|.
\end{align*}
\end{theorem}

Note that for $s \in \mathbb R$ one has $\|P(is)^{-1}\|=\|P(-is)^{-1}\|$ and, similarly, $\| R(is, \mathcal A_0)\|=\|R(-is, \mathcal A_0)\|$, so one may restrict attention  in the above estimates to large positive values of $s $.
Thus Theorem~\ref{resolvent_estimate} reduces the problem of estimating the decay rate of $(\mathcal T_0(t))_{t\ge0}$ to the stationary problem of obtaining norm estimates for $P(is)^{-1}$ for large $s>0$. The latter, however, is a very delicate matter and is intimately related to fine structure of $A$ and $B$. The problem of obtaining estimates for $P(is)^{-1}$ as $s\to\infty$  is far from being fully understood, even for comparatively simple models such as the damped wave equation discussed below. 

On the other hand,   resolvent estimates for $\mathcal A$, and hence  norm estimates for $P(\cdot)^{-1}$, can often be deduced from control-theoretic properties of the system governed by \eqref{equation} and/or joint spectral properties of $A$ and $B$. The first approach relies, in particular, on the study of various observability properties including Hautus-type tests, the second approach concerns so-called \emph{wave-packet conditions}.
Their advantages become apparent when dealing  with ``rough'' dampings, where there is less to be gained from employing microlocal techniques. Among the papers following this (often rather effective) ideology we mention, for example, \cite{BuZw04,RTTT05,Mi05,Mi12,AnLe14,JoLa20}.
To give a flavour of the relevant results, we formulate several resolvent estimates proved (in a more general context allowing also for unbounded operators $B$) in the recent paper \cite{CPSST19}. Similar (and sometimes stronger) results can be found in the papers mentioned above.

Let ${\boldsymbol A}$ be a skew-adjoint operator and ${\boldsymbol B}$ a bounded self-adjoint positive-definite operator on some Hilbert space. The pair $({\boldsymbol A}, {\boldsymbol B})$ is said to satisfy the \emph{non-uniform Hautus test} if there exist increasing functions $p$, $q: \mathbb R_+ \to (0,\infty)$ with $r_0>0$ such that
\begin{equation}\label{NUH}
  \| x\|^2\leq p(|s|) \| (is- {\boldsymbol A})x\|^2 + q(|s|) \| {\boldsymbol B} x\|^2 , \quad  x\in D({\boldsymbol A}),  s\in\mathbb R.
\end{equation}
Note that we may express our operator $\mathcal A$ in the form $\mathcal A = {\boldsymbol A} - {\boldsymbol B}$ where $${\boldsymbol A}=\left( \begin{array}{cc} 0 & I \\ -A & 0 \end{array} \right)\quad \mbox{and}\quad {\boldsymbol B}= \left( \begin{array}{cc} 0 & 0 \\ 0 & B \end{array} \right)$$ are of the requisite type. In this situation, the non-uniform Hautus test \eqref{NUH} can be rewritten (without essential loss) in terms of $A$ and $B$. Indeed, by \cite[Proposition 3.10]{CPSST19}, if $\tilde{p}$, $\tilde{q}:  \mathbb R_+\to (0,\infty)$ are increasing functions such that 
 \begin{equation} \label{NUH_WE}
   \|x\|^2\leq \tilde{p} (s) \| (s-A^{1/2})x\|^2 + \tilde{q} (s)\| B^{1/2} x\|^2,  \quad x \in D(A^{1/2}), s\geq 0,
 \end{equation}
 then $({\boldsymbol A},{\boldsymbol B})$ satisfies \eqref{NUH} with $q(s)=4 \tilde{q} (s)$ and $p(s)=C(\tilde{p}(s)+\tilde{q}(s))$ for some $C>0$. Conversely, if~\eqref{NUH} holds then~\eqref{NUH_WE} is satisfied for $\tilde{p}(s)=p(s)$ and $\tilde{q}(s)=q(s)/2$. This allows us to obtain the following resolvent estimate; see \cite[Proposition 3.10]{CPSST19}.
 
\begin{theorem}\label{res_observab}
If the pair $(A, B)$ satisfies \eqref{NUH_WE}, then  $i\mathbb R\subseteq \rho(\mathcal A)$ and $\| R(is, \mathcal A)\| \le C(\tilde{p}(|s|)+\tilde{q} (|s|))$ for all $s \in \mathbb R$ and some $C>0$.
\end{theorem}

Let $A$ be self-adjoint and let $\delta:\mathbb R\to (0,\infty)$ be a bounded function. Define  the $(s,\delta(s))$-\emph{wave-packet}  $\mathrm{WP}_{s,\delta(s)}(\mathcal A)$ as the spectral subspace of $A$ associated with the segment
$(s-\delta(s), s +\delta(s))$. 
The next statement provides an estimate for $\|R(is, \mathcal A)\|$ (which is proved via an intermediate estimate for $\|(-s^2 + is B + A)^{-1}\|$) in terms of \emph{joint} spectral properties of $A$ and $B$; see \cite[Theorem~3.8]{CPSST19}.

\begin{theorem}\label{WP}
If $\gamma: [0,\infty) \to (0,\infty)$ is a bounded function such that
\begin{equation}\label{WPa}
\|B^{1/2} x\|\ge \gamma (s) \|x\|, \qquad x \in \mathrm{WP}_{s, \delta(s)}(A^{1/2}),
\end{equation}
for all $s \ge 0$, then $i\mathbb R \subseteq \rho(\mathcal A)$ and $\|R(is, \mathcal A)\|\le C\gamma(|s|)^{-2}\delta(|s|)^{-2}$
for all $s \in \mathbb R$ and some $C>0$.
\end{theorem}

The wave-packet estimate \eqref{WPa}  can in turn be used to pass from a slightly different Hautus-type condition for $A$ and $B$ to a resolvent estimate for $\mathcal A$; see \cite[Proposition 3.9]{CPSST19}.

\begin{theorem} \label{Relation}
  If $p_S : \mathbb R_+\to(0,\infty)$ and $q_S: \mathbb R_+\to [r_0,\infty)$ with $r_0>0$ and $\eta:=\inf_{s\geq 0} p_S (s)(1+s)^2>0$ are such that
  \begin{equation}\label{Schh}
    \| x\|^2 \leq  p_S(s) \| (s^2-A)x \|^2 +  q_S(s) \| B^{1/2} x\|^2, \quad  x\in D(A),  s\geq 0,
  \end{equation}
  then  $i\mathbb R\subseteq \rho(\mathcal A)$ and $\|R(is,\mathcal A)\|\le C (1+s^2)  p_S (|s|)  q_S (|s|)$ for all $s\in\mathbb R$ and some $C>0$. 
\end{theorem}

Indeed, one may prove that if \eqref{Schh} holds then the conditions of Theorem~\ref{WP} are satisfied for $\delta(s) = c_0(2  p_S(|s|)(1+s^2))^{-1/2}$ and $\gamma(s) =(2  q_S(|s|))^{-1/2}$, with $c_0=\min(\sqrt{\eta},1/2)$.

Yet another ``non-uniform'' observability condition is given by the estimate
\begin{equation}\label{semigroup_ob}
C_\tau \|(1- {\boldsymbol A})^{-\beta}x\|^2 \le \int_{0}^{\tau} \|{\boldsymbol B} T_{\boldsymbol A}(t)\|^2 \, dt, \qquad x \in D({\boldsymbol A}), 
\end{equation}
for some fixed non-negative constants $\beta$, $\tau$, and $C_\tau$, where $(T_{\boldsymbol A}(t))_{t \in \mathbb R}$ is the unitary $C_0$-group generated by ${\boldsymbol A}$. This condition, which is used widely and often in a context which is similar to ours (see e.g. \cite{Du07,Mi12}), also leads to resolvent estimates, for example by relating it to  wave-packets conditions  of the form \eqref{WPa}.
For more on the asymptotic consequences of \eqref{semigroup_ob}, including polynomial (and more general) rates of energy decay for \eqref{equation}, we refer to \cite[Section 4]{CPSST19} and the references therein.

\section{Concluding remarks}

The aim of this survey has been to present recently developed abstract methods for the study of rates of decay of semigroup orbits. Unfortunately, we are unable in this short text to give a full account of the many exciting applications of semigroup methods and  the theory of resolvent estimates. In closing, however, we turn briefly to one particularly rich example, namely the damped wave equation
\begin{equation}\label{DWE}
 \begin{split}
  u_{tt} + b(x)u_t -\Delta u = 0 & \text{ in } (0,\infty ) \times \mathcal M, \quad u = 0 \, \text{ on } (0,\infty ) \times\partial \mathcal M , \\
  u(0,\cdot ) = u_0 & \text{ in } \mathcal M , \quad
  u_t (0,\cdot ) = u_1 \, \text{ in } \mathcal M ,
 \end{split}
\end{equation}
on a compact, connected, smooth Riemannian manifold $(\mathcal M,g)$ with (smooth) boundary $\partial M$ (which might be empty). Here $b\in L^\infty (\mathcal M)$ is a nonnegative function such that $b >0$ on a subset of $\mathcal M$ with positive measure,  and $\Delta=\Delta_g$ is the Laplace-Beltrami operator on $\mathcal M$. Of course, the damped wave equation is an example of a partial differential equation which can be rewritten as an abstract second order problem of the form \eqref{equation} in $X=L^2 (\mathcal M)$. In this case $A=-\Delta$ is the (negative) Laplace-Beltrami operator with Dirichlet boundary conditions on $L^2 (\mathcal M)$ and $B=b$ is a multiplication operator. Both are self-adjoint and positive semi-definite. Since $\mathcal M$ is connected, and by Calder\'on's unique continuation principle, non-zero eigenfunctions of $\Delta$ cannot vanish on an open set. Thus, if $b >0$ almost everywhere on an open set, then $bu \not= 0$ for every non-zero eigenfunction of $\Delta$. It follows that the operator matrix $$\mathcal A = \left( \begin{array}{cc} 0 & I \\ \Delta & -b \end{array} \right)$$ on the space $\mathcal X = H^1_0 (\mathcal M) \times L^2 (\mathcal M)$ (or more precisely, its restriction $\mathcal A_0$ to an appropriate subspace $\mathcal X_0$) has no spectrum on the imaginary axis and that $\mathcal A_0$ generates a semi-uniformly stable $C_0$-semigroup $(\mathcal T_0(t))_{t\ge0}$. Note that $\mathcal A = \mathcal A_0$ if and only if the boundary of the manifold $\mathcal M$ is nonempty.

For a solution $u$ of the damped wave equation we define the physical energy (the sum of the potential and the kinetic energy)  by
\begin{align*}
 E(u, t) & = \frac12 \int_{\mathcal M} |\nabla u|^2 + \frac12 \int_{\mathcal M} |u_t|^2 , 
\end{align*}
which corresponds precisely to the energy defined in the previous section. In particular, the decay rate of the energy $E(u,t)$ of classical solutions is the same as the decay rate of $\| \mathcal T_0 (t) \mathcal A_0^{-1}\|^2$. 
 
By using Carleman estimates (a somewhat quantified version of Calder\'on's unique continuation principle) one may show that the resolvent of $\mathcal A$ grows at most exponentially along the imaginary axis \cite[Th\'eor\`eme 1 (i)]{Le96}, and therefore the energy decays at least logarithmically \cite{Bq98} (see also Theorem \ref{semigrouprates}). This exponential resolvent bound is optimal in general: take $\mathcal M$ a $2$-sphere (hence $\mathcal M$ has strictly positive curvature) and $b$ a damping localised near a pole; see \cite[Th\'eor\`eme 1 (ii)]{Le96}. Consequently, one cannot in general hope to obtain faster than logarithmic energy decay. Depending on the geometry of the manifold and/or the regularity and the localisation of the support of the damping function $b$ one may, however, obtain better decay rates in some cases. 

If $\mathcal M = \mathbb T^2$ is the $2$-torus and $b\in L^\infty (\mathcal M)$, then the energy decays with the rate $O(t^{-\frac12})$, and the decay improves to $O(t^{-(1-\varepsilon)})$ for certain $b\in W^{k,\infty} (\mathcal M)$, no matter how small the support of $b$; see Anantharaman and L\'eautaud \cite[Theorems 2.3 and 2.6]{AnLe14}. On the other hand, for simple characteristic functions $b$, the optimal rate of decay is $O(t^{-\frac23})$ \cite{AnLe14,Sta17}, no matter how large the support of $b$. The fact that the regularity of $b$ may play a more important role than its support is also nicely illustrated in \cite{LeLe17,DaKl20,Kl19,Kl20}.

If the manifold $\mathcal M$ has strictly negative curvature and no boundary, and if $\gamma$ is a hyperbolic closed geodesic in $\mathcal M$ and  the damping $b$ vanishes in an $\varepsilon$-neighbourhood of $\gamma$ but is positive everywhere else, then the energy decays semi-uniformly with the rate $O(e^{-\alpha \sqrt{t}})$ by \cite[Theorem 1]{BuCh15}; see \cite[Section 5]{CSVW14}, \cite[Corollary A2]{Ri14}, \cite{Sc10}, \cite[Theorem 1]{Sc11} and \cite[Theorem 1.1]{Jin20} for related results. The results in \cite{Le96}, \cite[Theorem 1]{Sc11} and \cite{Jin20} are particularly interesting since they show that the energy of classical solutions may decay exponentially without the semigroup $(\mathcal T_0(t))_{t\ge0}$ being exponentially stable; this provides another counterexample to the equality of the exponential growth bound and the spectral bound. 

Nevertheless, exponential stability of the semigroup $(\mathcal T_0(t))_{t\ge0}$ does occur, for instance in the simple situations where $\mathcal M$ is a compact interval \cite{CoZu94} or when $b\geq b_0$ for some constant $b_0 >0$, that is, the damping is active on the whole manifold. Building on the the principle of \emph{propagation of singularities} (see for instance H\"ormander \cite{Ho71}), Ralston \cite{Ra69} (see also \cite{Ra82}) showed that the energy of solutions to wave-type equations may localise along a single generalised geodesic (or, in different terminology, along a single ray of geometric optics). As a result, it turns out to be possible to characterise exponential stability of $(\mathcal T_0(t))_{t\ge0}$ in terms of the relationship between the support of $b\in C(\mathcal M)$ and the generalised geodesics or, more precisely, in terms of the so-called {\em geometric control condition}, which very roughly means that every generalised geodesic enters the set $\{ b>0\}$ in finite time; we do not go into details here. Sufficiency of the geometric control condition was proved for smooth dampings and manifolds without boundaries by Rauch and Taylor \cite{RaTa74}; see also \cite{Ra69,RaTa75}. As a culmination of a series of preceding works \cite{BaLeRa88}, \cite{Li88a} and \cite{BaLeRa89}, the result was extended in \cite{BaLeRa92} by Bardos, Lebeau and Rauch to the case of manifolds \emph{with} boundaries in the sense that they formulated a sufficient geometric control condition which is also close to being necessary. Later on the approach of \cite{BaLeRa92} was simplified and the result was generalised to continuous dampings in \cite{BuGe97} by means of control-theoretical arguments (in the formally less general setting of smooth domains).
Note also that a detailed treatment of geodesic flows from the point of view of fine observability estimates for wave equations may be found in the recent paper \cite{LaLe16}, which moreover improves on several similar arguments given in \cite{BaLeRa92}. 


Finally, we mention that rates of decay for the \emph{derivatives} of the energy
in the general setting of (5.1) were obtained by means of the abstract $L^p$-Tauberian Theorem \ref{x0} in \cite[Corollary 6.5]{BaBoTo16}.

Many of the results above rely on abstract semigroup theory, often on resolvent bounds. However, we emphasize that in special geometries and for special damping functions it is sometimes possible to derive so-called integral inequalities for the energy. These integral inequalities, which are typically based on sophisticated  \emph{multiplier methods}  and clever estimates of various integrals using monotonicity or convexity assumptions, lead directly to decay rate estimates for the associated energies without the passage through the resolvent estimates. There is a vast literature literature on this approach; we restrict ourselves to mentioning just a few illustrative examples, such as  J.-L. Lions \cite{Lio88}, Komornik \cite{Ko94}, Martinez \cite{Mr99,Mr99a} and Alabau \cite{Al05}. This approach has the considerable advantage of being applicable also to non-linear equations.  A comparatively recent and comprehensive survey of applications of multiplier methods to hyperbolic equations may be found in \cite{Al12}.

The classical damped wave equation discussed here is, of course, just one (prototypical) example of how abstract semigroup results may be applied to concrete PDEs. We conclude by mentioning  a small selection of other interesting applications of semigroup methods to differential equations, along with sample references:  damped wave equations on unbounded domains/manifolds \cite{AlKh02,Kh03, BuJo16,Da16,Nis16,JoRo18,MaRo18};  local energy decay for damped wave equations \cite{BoRo14,Ro18};  Klein-Gordon and Kelvin-Vogt type equations \cite{Wu17,AmHaRo20,Bu20}; energy decay for non-linear damped wave equations \cite{Bd05,Ph11,JoLa13,JoLa20}; vectorial damped wave equations \cite{Kl18}; damped wave equations with unbounded and/or indefinite dampings\cite{FrZu96,BeRa00,AbMeNi13,FrSiTr18,ArSi20}; viscoelastic boundary dampings \cite{Sta18}; wave equations with periodic (or even general non-stationary) dampings \cite{LLTT17,JoRo18}; fractional damped wave equations \cite{Gr19,Gr19a}; damped wave equations on manifolds with rough metrics \cite{Tay15}. Even though the subject of damped wave equations is already vast, we hope and expect that the stream of substantial advances in this area, whether obtained by abstract techniques or otherwise, will continue for many years to come.\\

{\bf Acknowledgements.} The authors kindly thank the two anonymous referees for their numerous helpful comments and suggestions.



 \def\cprime{$'$}
  \def\ocirc#1{\ifmmode\setbox0=\hbox{$#1$}\dimen0=\ht0 \advance\dimen0
  by1pt\rlap{\hbox to\wd0{\hss\raise\dimen0
  \hbox{\hskip.2em$\scriptscriptstyle\circ$}\hss}}#1\else {\accent"17 #1}\fi}
  \def\cprime{$'$} \def\cprime{$'$} \def\cprime{$'$}

\end{document}